\documentclass[12pt,oneside,reqno]{amsart}
\usepackage{amsmath}
\usepackage{graphicx}
\usepackage{mathrsfs}
\usepackage{stmaryrd}
\usepackage{amsfonts}
\usepackage{enumerate,amsmath,amssymb,amsthm}
\newcommand{\be}{\begin{eqnarray}}
\newcommand{\ee}{\end{eqnarray}}
\newcommand{\ce}{\begin{eqnarray*}}
\newcommand{\de}{\end{eqnarray*}}
\newtheorem{theorem}{Theorem}[section]
\newtheorem{lemma}[theorem]{Lemma}
\newtheorem{remark}[theorem]{Remark}
\newtheorem{definition}[theorem]{Definition}
\newtheorem{proposition}[theorem]{Proposition}
\newtheorem{Examples}[theorem]{Examples}
\newtheorem{corollary}[theorem]{Corollary}

\def\a{\alpha}

\def\b{\beta}
\def\p{\partial}
\def\d{\delta}
\def\eps{\epsilon}

\def\[{{\Big[}}
\def\]{{\Big]}}
\def\<{{\langle}}
\def\>{{\rangle}}
\def\({{\Big(}}
\def\){{\Big)}}

\def\bx{{\mathbf{x}}}

\def\dif{{\mathord{{\rm d}}}}

\def\b{\beta}

\def\no{\nonumber}
\def\bt{\begin{theorem}}
\def\et{\end{theorem}}
\def\bl{\begin{lemma}}
\def\el{\end{lemma}}
\def\br{\begin{remark}}
\def\er{\end{remark}}
\def\bx{\begin{Examples}}
\def\ex{\end{Examples}}
\def\bd{\begin{definition}}
\def\ed{\end{definition}}
\def\bp{\begin{proposition}}
\def\ep{\end{proposition}}
\def\bc{\begin{corollary}}
\def\ec{\end{corollary}}
\def\cA{{\mathcal A}}
\def\cB{{\mathcal B}}

\def\cD{{\mathcal D}}

\def\cF{{\mathcal F}}

\def\cO{{\mathcal O}}

\def\mB{{\mathbb B}}
\def\mC{{\mathbb C}}
\def\mD{{\mathbb D}}
\def\mE{{\mathbb E}}

\def\mN{{\mathbb N}}

\def\mP{{\mathbb P}}

\def\mR{{\mathbb R}}
\def\mS{{\mathbb S}}

\def\mU{{\mathbb U}}

\def\mZ{{\mathbb Z}}

\def\geq{\geqslant}
\def\leq{\leqslant}

\allowdisplaybreaks

\begin{document}
\title{Large Deviations for Multi-valued Stochastic Differential Equations$^*$}
\author{Jiagang Ren$^{1}$, Siyan Xu$^{2}$, Xicheng Zhang$^{2,3}$}


\date{}
\dedicatory{$^{1}$ School of Mathematics and Computational Science,
Sun Yat-Sen University,\\
 Guangzhou, Guangdong 510275, P.R.China\\
$^{2}$School of Mathematics and statistic,
Huazhong University of Science and Technology,\\
Wuhan, Hubei 430074, P.R.China\\
$^{3}$School of Mathematics and Statistics,\\
The University of New South Wales, Sydney, 2052, Australia\\
Emails: J. Ren: renjg@mail.sysu.edu.cn\\
S. Xu: xsy\_00@hotmail.com\\
 X.Zhang: xichengzhang@gmail.com}


\thanks{This work is supported by NSF of China (No. 10871215). The third author
is also supported by ARC Discovery grant DP0663153 of Australia.}
\begin{abstract}
We prove a large deviation principle of Freidlin-Wentzell's type for the
multivalued stochastic differential equations with monotone drifts,
which in particular contains a class of SDEs with reflection in a convex domain.
\end{abstract}

\maketitle

\section{Introduction}

Consider the following  multivalued stochastic differential equation
(MSDE in short):
\be\label{eq0} \left\{
\begin{array}{ll}
\dif X(t)\in b(X(t))\dif t+\sigma(X(t))\dif W(t)-A(X(t))\dif t,\\
X(0)=x\in \overline{D(A)},
\end{array}
\right.
\ee
where $A$ is a multivalued maximal monotone operator, which will be described below,
$W(t)=\{W^k(t), t\geq 0, k\in\mN\}$ is a sequence of independent standard Brownian motions
on a filtered probability space $(\Omega,\cF,P; (\cF_t)_{t\geq 0})$,
$b: \mR^m\to\mR^m$ and $\sigma:\mR^m\to\mR^m\times l^2$ are two continuous functions,
$l^2$ stands for the Hilbert space of square summable sequences of real numbers.

This type of MSDE was first studied by C\'epa in \cite{ce,ce2}.
He proved that if $b$ and $\sigma$ are Lipschitz continuous,
then there exists a unique pair of processes
$(X(t),K(t))$ such that
$$
X(t)=x+\int^t_0b(X(s))\dif s+\int^t_0\sigma(X(s))\dif W(s)-K(t),
$$
where $K(t)$ is a process of finite variation (see Definition
\ref{defsolution} below for more details). Recently, Zhang
\cite{zhang} extended C\'epa's result to the infinite dimensional
case, and relaxed the Lipschitz assumption on $b$ to the monotone
case. It should be noted that when $A$ is the subdifferential of the
indicator function of a convex subset of $\mR^m$, the above MSDE is
the same as the usual SDE with reflecting  boundary in a convex
domain  (cf. \cite{An, wa}). Moreover, since the subdifferential of
any lower semicontinuous convex function is a maximal monotone
operator, C\'epa's result can also be used to deal with the SDE with
discontinuous coefficients. It is well known that there are many
literatures to investigate the SDEs with reflecting boundary since
the solutions of a large class of PDEs with Neumann boundary and
mixed boundary conditions can be represented by the solution of such
SDEs (cf. \cite{An}).

We now consider the following small perturbation of Eq.(\ref{eq0}):
\be \left\{
\begin{array}{ll}
\dif X^\eps(t)\in b(X^\eps(t))\dif t+\sqrt{\eps}
\sigma(X^\eps(t))\dif W(t)-A(X^\eps(t))\dif t,\\
X^\eps(0)=x\in \overline{D(A)},\ \ \eps\in(0,1].
\end{array}
\right. \label{pmsde}
\ee
The solution of this equation is denoted by $(X^\eps(t,x),K^\eps(t,x))$.
We want to establish the large deviation principle of the law of
$X^\eps(t,x)$ in the space $\mS:=C([0,T]\times\overline{D(A)};\overline{D(A)})$, namely,
the asymptotic estimates of probabilities $P(X^\eps\in \Gamma)$,
where $\Gamma\in \cB(\mS)$.

In \cite{An}, Anderson and Orey considered the same small random
perturbation for the dynamical system with reflecting boundary in
smooth domain, and obtained the Freidlin-Wentzell's large deviation
estimates in $C([0,T];\overline{D(A)})$. They assumed that the
coefficients are bounded and Lipschitz continuous, and the diffusion
coefficient is non-degenerate. Using the contraction principle,
C\'epa \cite{ce} only  considered  the large deviation principle of
one dimensional case based on an explicit construction of the
solution (cf.\cite{xu}). The multi-dimensional case is still open.
Compared with the usual SDE, i.e., $A=0$, most of the difficulties
come from the presence of the process of finite variation,  $K(t)$.
One only knows that $t\mapsto K(t)$ is continuous, and could not
prove any further regularity such as H\"older continuity. Therefore,
the classical method of time discretization is almost inapplicable
(cf. \cite{DZ}).

Our method is based on the recently well developed weak convergence
approach due to Dupuis and Ellis \cite{Du-El} (see also \cite{bd,
bd0}). This method has been proved to be very effective for various
systems (cf. \cite[etc.]{Rz2, Zh2, Bu-Du-Ma, Zh0, Liu,Rz3}). In the
situation considered in the present paper, however, since we cannot
prove the following uniform estimate as in \cite{Rz2}: for any
$p\geq 2$ and $s,t\in[0,T]$, $x,y\in \overline{D(A)}$
$$
\sup_{\eps\in(0,1)}\mE|X^\eps(t,x)-X^\eps(s,y)|^{2p}\leq C(|t-s|^p+|x-y|^{2p}),
$$
we cannot obtain the tightness of the laws of $X^\eps(t,x)$ in
$\mS$. Some technical difficulties for verifying the conditions {\bf
(LD)$_1$} and {\bf (LD)$_2$} below need to be overcome.

In Section 2, we recall some well known facts about the MSDE and a criterion for Laplace principle.
In Section 3, we present our main result and give a detailed proof.
Throughout the paper, $C$ with or without indexes will denote different constants
(depending on the indexes) whose values are not important.

\section{Preliminaries}

We first give some notions and notations about multivalued operators.
Let $2^{\mR^m}$ be the set of all subsets of $\mR^m$. A map $A:\mR^m\to 2^{\mR^m}$
is called a multivalued operator. Given such a multivalued operator $A$, define:
\ce
D(A)&:=&\{x\in \mR^m:A(x)\neq \emptyset\}, \\
\mathrm{Im}(A)&:=&\cup_{x\in D(A)}A(x),\\
\mathrm{Gr}(A)&:=&\{(x,y)\in \mR^{2m}: x\in \mR^m, y\in A(x)\}.
\de

We recall the following definitions.
\bd\label{Def1}
(1) A multivalued operator $A$ is called monotone if
\ce
\<y_1-y_2,x_1-x_2\>_{\mR^m}\geq 0, \quad \forall (x_1,y_1),(x_2,y_2)\in \mathrm{Gr}(A).
\de
(2) A monotone operator A is called maximal monotone if for each
$(x,y)\in \mathrm{Gr}(A)$,
$$
\<y-y',x-x'\>_{\mR^m} \geq 0,\quad\forall(x',y')\in \mathrm{Gr}(A).
$$
\ed

\bx
Suppose that $\cO$ is a closed convex subset of $\mR^m$, and
$I_\cO$ is the indicator function of $\cO$, i.e,
$$
I_\cO(x)=\left\{
\begin{array}{ll}
0, &~\mbox{if } \ x \in \cO,\\
+\infty,&~\mbox{if }\ ~x \notin \cO.
\end{array}
\right.
$$
The subdifferential of $I_\cO$ is given by
\ce
\p I_\cO(x)&£º=&\{y\in \mR^m:\<y,x-z\>_{\mR^m}\geq 0, \forall z\in \cO\}\\
&=&\left\{
\begin{array}{ll}
\emptyset, ~&\mbox{if }\ ~x \notin \cO,\\
\{0\},~&\mbox{if }\ ~x \in \mathrm{Int}(\cO),\\
\Lambda_x, ~&\mbox{if }\ ~x\in \p \cO,
\end{array}
\right.
\de
where $\mathrm{Int}(\cO)$ is the interior of $\cO$ and $\Lambda_x$ is the exterior normal cone at $x$.
One can check that $\p I_\cO$ is a multivalued maximal monotone operator in the sense of Definition
\ref{Def1}.
\ex

We now give the precise definition of the solution to Eq.(\ref{eq0}).
\bd\label{defsolution}
A pair of continuous and ($\cF_t$)-adapted processes $(X,K)$ is called
a solution of Eq.(1) if
\begin{enumerate}[(i)]
\item $X(0)=x$, and for all $t\geq 0$, $X(t)\in \overline{D(A)}\quad a.s.$;

\item $K(0)=0$ a.s. and $K$ is of finite variation;

\item $\dif X(t)=b(X(t))\dif t+\sigma(X(t))\dif W(t)-\dif K(t)$, $0\leq t<\infty, \quad a.s.$;

\item for any continuous and ($\cF_t$)-adapted processes $(\a,\b)$ with
$$
(\a(t),\b(t))\in \mathrm{Gr}(A),\quad \forall t\in[0,+\infty),
$$
the measure
$$
\<X(t)-\a(t), \dif K(t)-\b(t)\dif t\>\geq 0 \quad a.s..
$$
\end{enumerate}
\ed

We now recall an abstract criterion for Laplace principle, which is
equivalent to the large deviation principle (cf.
\cite{bd0,Bu-Du-Ma,Zh2}). It is well known that there exists a
Hilbert space so that $l^2\subset\mU$ is Hilbert-Schmidt with
embedding operator $J$ and $\{W^k(t),k\in\mN\}$ is a  Brownian
motion with values in $\mU$, whose covariance operator is given by
$Q=J\circ J^*$. For example, one can take $\mU$ as the completion of
$l^2$ with respect to the norm generated by the scalar product
$$
\<h,h'\>_\mU:=\left(\sum_{k=1}^\infty \frac{h_k h'_k}{k^2}\right)^{\frac{1}{2}}, \ \ h,h'\in l^2.
$$

For a Polish space $\mB$, we denote by $\cB(\mB)$ its   Borel
$\sigma$-field, and by $\mC_T(\mB)$ the continuous function space
from $[0,T]$ to $\mB$, which is endowed with the uniform distance so
that $\mC_T(\mB)$ is still a Polish space. Define \be
\ell^2_T:=\left\{h=\int^\cdot_0\dot h(s)\dif s: ~~\dot h\in
L^2(0,T;l^2)\right\}\label{H2} \ee with the norm
$$
\|h\|_{\ell^2_T}:=\left(\int^T_0\|\dot h(s)\|_{l^2}^2\dif s\right)^{1/2},
$$
where the dot denotes the generalized derivative.
Let $\mu$ be the law of the Brownian motion $W$ in $\mC_T(\mU)$. Then
$$
(\ell^2_T,\mC_T(\mU),\mu)
$$
forms an abstract Wiener space.

For $T,N>0$,  set
$$
\cD_N:=\{h\in \ell^2_T: \|h\|_{\ell^2_T}\leq N\}
$$
and
\be
\cA^T_N:=\left\{
\begin{aligned}
&\mbox{ $h: [0,T]\to l^2$ is a continuous and
$(\cF_t)$-adapted }\\
&\mbox{ process, and for almost all $\omega$},\ \ h(\cdot,\omega)\in\cD_N
\end{aligned}
\right\}.\label{Op2}
\ee
We equip $\cD_N$ with the  weak convergence topology in
$\ell^2_T$. Then
\be
\mbox{$\cD_N$ is metrizable as a compact Polish space}.\label{Metr}
\ee

Let $\mS$ be a Polish space. A function $I: \mS\to[0,\infty]$ is given.
\bd
The function $I$  is called a rate function if for every $a<\infty$, the set
$\{f\in\mS: I(f)\leq a\}$ is compact in $\mS$.
\ed

Let $\{Z^\eps: \mC_T(\mU)\to\mS,\eps\in(0,1)\}$
be a family of measurable mappings. Assume that
there is a measurable map $Z_0: \ell^2_T\mapsto \mS$ such that
\begin{enumerate}[{\bf (LD)$_\mathbf{1}$}]
\item
For any $N>0$, if a family $\{h_\eps, \eps\in(0,1)\}\subset\cA^T_N$ (as random variables in $\cD_N$)
converges in distribution to $h\in \cA^T_N$, then
for some subsequence $\eps_k$, $Z^{\eps_k}\Big(\cdot+\frac{h_{\eps_k}(\cdot)}{\sqrt{\eps_k}}\Big)$
converges in distribution to $Z_0(h)$ in $\mS$.
\end{enumerate}
\begin{enumerate}[{\bf (LD)$_\mathbf{2}$}]
\item
 For any $N>0$, if $\{h_n,n\in\mN\}\subset \cD_N$ weakly converges to $h\in\ell^2_T$,
then for some subsequence $h_{n_k}$, $Z_0(h_{n_k})$ converges to $Z_0(h)$ in $\mS$.
\end{enumerate}

For each $f\in\mS$, define
\be
I(f):=\frac{1}{2}\inf_{\{h\in\ell^2_T:~f=Z_0(h)\}}\|h\|^2_{\ell^2_T},\label{ra}
\ee
where $\inf\emptyset=\infty$ by convention. Then under {\bf (LD)$_\mathbf{2}$},
$I(f)$ is a  rate function.

We recall the following result due to \cite{bd0} (see also \cite[Theorem 4.4]{Zh1}).
\bt\label{Th2}
Under {\bf (LD)$_\mathbf{1}$} and {\bf (LD)$_\mathbf{2}$},
$\{Z^\eps,\eps\in(0,1)\}$ satisfies the Laplace principle with
the rate function $I(f)$ given by (\ref{ra}). More precisely, for each real bounded continuous
function  $g$ on $\mS$:
\be
\lim_{\eps\rightarrow 0}\eps\log\mE^{\mu}\left(\exp\left[-\frac{g(Z^\eps)}{\eps}\right]\right)
=-\inf_{f\in\mS}\{g(f)+I(f)\}.\label{La}
\ee
In particular, the family of $\{Z^\eps,\eps\in(0,1)\}$
satisfies  the large deviation principle in $(\mS,\cB(\mS))$ with the rate function $I(f)$.
More precisely, let $\nu_\eps$ be the law of $Z^\eps$ in $(\mS,\cB(\mS))$,
then for any $B \in\cB(\mS)$
$$
-\inf_{f\in B^o}I(f)\leq\liminf_{\eps\rightarrow 0}\eps\log\nu_\eps(B)
\leq\limsup_{\eps\rightarrow 0}\eps\log\nu_\eps(B)\leq -\inf_{f\in \bar B}I(f),
$$
where the closure and the interior are taken in $\mS$,
and $I(f)$ is defined by (\ref{ra}).
\et

\section{Main Result and Proof}

We assume that
\begin{enumerate}[{\bf (H1)}]
\item $A$ is a maximal monotone operator with non-empty interior, i.e., Int$(D(A))\not=\emptyset$;
\item  $\sigma$ and $b$ are continuous functions
and satisfy that for some $C_\sigma, C_b>0$ and all $x,y\in\mR^m$
\ce
\|\sigma(x)-\sigma(y)\|_{L_2(l^2;\mR^m)}&\leq& C_\sigma|x-y|,\\
\<x-y,b(x)-b(y)\>_{\mR^m}&\leq& C_b|x-y|^2,
\de
where $L_2(l^2;\mR^m)$ denotes the Hilbert-Schmidt space and $|\cdot|$ denotes the norm in $\mR^m$,
and for some $C_b'>0$ and $n\in\mN$
$$
|b(x)|\leq C_b'(1+|x|^n).
$$
\end{enumerate}

It is well known that under {\bf (H1)} and {\bf (H2)}, there exists a unique solution
$(X^\eps,K^\eps)$ to Eq.(\ref{pmsde}) in the sense of Definition \ref{defsolution} (cf. \cite{zhang}).
Our main result is stated as follows:
\bt\label{Main}
Assume that {\bf (H1)} and {\bf (H2)} hold.
Then the family of $\{X^\eps(t,x),\ \eps\in(0,1)\}$ satisfies the large deviation principle
in $\mS:=C([0,T]\times \overline{D(A)};\overline{D(A)})$ with the rate function given by
\be
I(f):=\frac{1}{2}\inf_{\{h\in\ell^2_T:~f=X^h\}}\|h\|^2_{\ell^2_T},\label{rae}
\ee
where $X^h(t,x)$ solves the following equation:
\ce
\dif X^h(t)\in b(X^h(t))\dif t+
\sigma(X^h(t))\dot h(t)\dif t-A(X^h(t))\dif t,\ \ X^h(0)=x.
\de
\et
For proving this result, by Theorem \ref{Th2},
the main task is to verify {\bf (LD)$_\mathbf{1}$}
and {\bf (LD)$_\mathbf{2}$} with
$$
\mS:=C([0,T]\times \overline{D(A)};\overline{D(A)}),\ \ Z^\eps=X^\eps,\ \ Z_0(h)=X^h.
$$
This will be done in Lemmas \ref{Le6} and \ref{Le7} below.

\vspace{5mm}

Let $h_\eps\in\cA^T_N$ converge almost surely to $h\in\cA^T_N$
as random variables in $\ell^2_T$, and $(X^{\eps,h_\eps}, K^{\eps,h_\eps})$
solve the following control equation:
\be
 X^{\eps,h_\eps}(t)&=& x+\int^t_0b(X^{\eps,h_\eps}(s))\dif s+
 \int^t_0\sigma(X^{\eps,h_\eps}(s))\dot{h}_\eps(s)\dif s\no\\
&&+\sqrt{\eps}\int^t_0\sigma(X^{\eps,h_\eps}(s))\dif W(s)-K^{\eps,h_\eps}(t),\label{NS0}
\ee
which can be solved by Girsanov's theorem,
and $(X^h, K^h)$ solve the following deterministic equation:
\be
 X^h(t)= x+\int^t_0b(X^h(s))\dif s
+\int^t_0\sigma(X^h(s))\dot{h}(s)\dif s-K^h(t).\label{Es2}
\ee

Let $|K|_t^s$ denote the total variation of $K$ on $[s,t]$.
We recall the following result due to C\'epa \cite{ce2}
(see also \cite[Propositions 3.3 and 3.4]{zhang}).
\bp\label{bp}
Under {\bf (H1)}, there exist $a\in \mR^m, \gamma>0, \mu \geq 0$ such that
for any pair of $(X,K)$ with the property (iv) of Definition \ref{defsolution}
and all $0\leq s<t\leq T$
\be
\int_s^t\<X(r)-a,\dif K(r)\>_{\mR^m}\geq \gamma|K|_t^s
-\mu\int_s^t|X(r)-a|\dif r-\gamma \mu(t-s).\label{Es7}
\ee
Moreover, for any pairs of $(X,K)$ and
$(\tilde X,\tilde K)$ with the property (iv) of Definition \ref{defsolution}
\be
\<X(t)-\tilde X(t),\dif K(t)-\dif\tilde K(t)\>_{\mR^m}\geq 0.\label{Es8}
\ee
\ep
Using this property, we first prove the following uniform estimates.
\bl\label{Le5}
For any $p\geq 1$, there exists $C_{p,T,N}>0$ such that for any $\eps\in(0,1)$
and $x,y\in\overline{D(A)}$
\be
\mE\left(\sup_{t\in[0,T]}|X^{\eps,h_\eps}(t,x)- X^{\eps,h_\eps}(t,y)|^{2p}\right)
\leq C_{p,T,N}|x-y|^{2p}.\label{Es1}
\ee
\el
\begin{proof}
Set
$$
Z_\eps(t):=X^{\eps,h_\eps}(t,x)- X^{\eps,h_\eps}(t,y)
$$
and
$$
\Lambda(s):=\sigma(X^{\eps,h_\eps}(s,x))-\sigma(X^{\eps,h_\eps}(s,y)).
$$
By It\^o's formula, {\bf (H2)} and (\ref{Es8}), we have for any $p\geq 1$
\ce
|Z_\eps(t)|^{2p}
&=&|Z_\eps(0)|^{2p}+2p\int_0^t|Z_\eps(s)|^{2p-2}
\<Z_\eps(s),b(X^{\eps,h_\eps}(s,x))- b(X^{\eps,h_\eps}(s,y))\>_{\mR^m}\dif s\\
&&+2p\int_0^t|Z_\eps(s)|^{2p-2}\<Z_\eps(s),\Lambda(s)\dot{h}_\eps(s)\>_{\mR^m}\dif s\\
&&+2p\sqrt{\eps}\int_0^t|Z_\eps(s)|^{2p-2}
\<Z_\eps(s),\Lambda(s)\dif W(s)\>_{\mR^m}\\
&&-2p\int_0^t|Z_\eps(s)|^{2p-2}\<Z_\eps(s),\dif K^{\eps,h_\eps}(s,x)-
\dif K^{\eps,h_\eps}(s,y)\>_{\mR^m}\\
&&-p\int_0^t|Z_\eps(s)|^{2p-2}(\|\Lambda(s)\|^2
+2(p-1)\<Z_\eps(s),\Lambda(s)\Lambda^*(s)Z_\eps(s)\>_{\mR^m}/|Z_\eps(s)|^2)\dif s\\
&\leq&|Z_\eps(0)|^{2p}+C\int_0^t|Z_\eps(s)|^{2p}\dif s
+2p\int_0^t|Z_\eps(s)|^{2p}\cdot\|\dot{h}_\eps(s)\|_{l^2}\dif s\\
&&+2p\sqrt{\eps}\int_0^t|Z_\eps(s)|^{2p-2}
\<Z_\eps(s),\Lambda(s)\dif W(s)\>_{\mR^m}.
\de

Set
$$
g(t):=\mE\left(\sup_{s\in[0,t]}|Z_\eps(s)|^{2p}\right).
$$
By BDG's inequality and Young's inequality,
we have for any $\d>0$
\be
&&\mE\left|\sup_{t'\in[0,t]}\int_0^{t'}|Z_\eps(s)|^{2p-2}
\<Z_\eps(s),\Lambda(s)\dif W(s)\>_{\mR^m}\right|\no\\
&&\qquad\qquad\leq C\mE\left(\int_0^{t}|Z_\eps(s)|^{4p-4}
\|\Lambda(s)^*Z_\eps(s)\|_{l^2}^2\dif s\right)^{1/2}\no\\
&&\qquad\qquad\leq C\mE\left(\sup_{s\in[0,t]}|Z_\eps(s)|^{2p}\int_0^{t}|Z_\eps(s)|^{2p}
\dif s\right)^{1/2}\no\\
&&\qquad\qquad\leq \d\cdot g(t)+C_\d\int_0^{t}\mE|Z_\eps(s)|^{2p}\dif s.\label{Es3}
\ee
Similarly, we have
\be
\mE\left|\int_0^t|Z_\eps(s)|^{2p}\cdot\|\dot{h}_\eps(s)\|_{l^2}\dif s\right|
&\leq& N\mE\left(\int_0^t|Z_\eps(s)|^{4p}\dif s\right)^{1/2}\no\\
&\leq&\d\cdot g(t)+C_{\d,N}\int_0^{t}\mE|Z_\eps(s)|^{2p}\dif s.\label{Es4}
\ee
Letting $\d=1/4$ in (\ref{Es3}) and (\ref{Es4}) and combining the above calculations,
we get
$$
g(t)\leq g(0)+\frac{1}{2}g(t)+C\int^t_0\mE|Z_\eps(s)|^{2p}\dif s
\leq 2g(0)+2C\int^t_0g(s)\dif s,
$$
which gives the desired estimate by Gronwall's inequality.
\end{proof}

\bl
For any $p\geq 1$ and $x\in\overline{D(A)}$, there exists $C_{p,T,N,x}>0$
such that  for any $\eps\in[0,1)$
\be
\mE\left(\sup_{t\in[0,T]}|X^{\eps,h_\eps}(t,x)|^{2p}\right)
+\mE|K^{\eps,h_\eps}(\cdot,x)|^0_T\leq C_{p,T,N,x}.\label{Es5}
\ee
\el
\begin{proof}
First of all, as in the proof of Lemma \ref{Le5} we can prove that
\be \mE\left(\sup_{t\in[0,T]}|X^{\eps,h_\eps}(t,x)|^{2p}\right)\leq
C_{p,T,N,x}.\label{Es10} \ee Let $a\in\mathrm{Int}(D(A))$ be as in
Proposition \ref{bp}. By It\^o's formula, (\ref{Es7}) and {\bf
(H2)}, we have \ce
\frac{1}{2}|X^{\eps,h_\eps}(t)-a|^{2}&=&\frac{1}{2}|x-a|^{2}
+\int_0^t\<X^{\eps,h_\eps}(s)-a,b(X^{\eps,h_\eps}(s))\>_{\mR^m}\dif s\\
&&+\int_0^t\<X^{\eps,h_\eps}(s)-a,\sigma(X^{\eps,h_\eps}(s))\dot{h}_\eps(u)\>_{\mR^m}\dif s\\
&&+\sqrt{\eps}\int_0^t\<X^{\eps,h_\eps}(s)-a,\sigma(X^{\eps,h_\eps}(s))\dif W(s)\>_{\mR^m}\\
&&-\int_0^t\<X^{\eps,h_\eps}(s)-a,\dif K^{\eps,h_\eps}(s)\>_{\mR^m}\\
&&+\frac{\eps}{2}\int_0^t\|\sigma(X^{\eps,h_\eps}(s))\|^2_{L_2(l^2;\mR^m)}\dif s\\
&\leq&\frac{1}{2}|x-a|^{2}+C_b\int^t_0|X^{\eps,h_\eps}(s)-a|^2\dif s\\
&&+\int_0^t\<X^{\eps,h_\eps}(s)-a, b(a)\>_{\mR^m}\dif s\no\\
&&+N\left(\int_0^t|\sigma(X^{\eps,h_\eps}(s))^*(X^{\eps,h_\eps}(s)-a)|^2\dif s\right)^{1/2}\no\\
&&+\sqrt{\eps}\int_0^t\<X^{\eps,h_\eps}(s)-a,\sigma(X^{\eps,h_\eps}(s))\dif W(u)\>_{\mR^m}
+\mu\gamma t\no\\
&&-\gamma|K^{\eps,h_\eps}|_t^0+\mu\int_0^t|X^{\eps,h_\eps}(s)-a|\dif s\\
&&+\frac{C^2_\sigma\eps}{2}\int_0^t(|X^{\eps,h_\eps}(s)|+\sigma(0))^2\dif s.
\de
The desired estimate now follows by (\ref{Es10}).
\end{proof}

Define
\be
w_\eps(t,x):=\int^t_0 \sigma(X^h(s,x))(\dot h_\eps(s)-\dot h(s))\dif s.\label{Def}
\ee
The following lemma is easy by Ascoli-Arzela's lemma.
\bl\label{Le3}
$w_\eps(\cdot,x)$ converges a.s. to zero in $C([0,T],\overline{D(A)})$.
\el

We now prove the following key lemma.
\bl\label{Le4}
$ X^{\eps,h_\eps}$ defined by (\ref{NS0}) converges in probability to $ X^h$ defined by
(\ref{Es2}) in $\mS$.
\el
\begin{proof} Set $v_\eps(t):=v_\eps(t,x):= X^{\eps,h_\eps}(t,x)- X^h(t,x)$. Then
\ce
 v_\eps(t)&=&K^{\eps,h_\eps}(t)-K^h(t)+\int^t_0(b(X^{\eps,h_\eps}(s))-b(X^h(s)))\dif s\no\\
&&+\int^t_0(\sigma(X^{\eps,h_\eps}(s))\dot h_\eps(s)-\sigma(X^h(s))\dot h(s))\dif s\\
&&+\sqrt{\eps}\int^t_0\sigma(X^{\eps,h_\eps}(s))\dif W(s).
\de
By It\^o's formula, we have
\ce
| v_\eps(t)|^2&=&2\int^t_0\<v_\eps(s),\dif K^{\eps,h_\eps}(s)-\dif K^h(s)\>_{\mR^m}\\
&&+2\int^t_0\< v_\eps(s),b(X^{\eps,h_\eps}(s))-b(X^h(s))\>_{\mR^m}\dif s\\
&&+2\int^t_0\<v_\eps(s),(\sigma(X^{\eps,h_\eps}(s))-\sigma(X^h(s)))\dot h_\eps(s)\>_{\mR^m}\dif s\\
&&+2\int^t_0\<v_\eps(s),\sigma(X^h(s))(\dot h_\eps(s)-\dot h(s))\>_{\mR^m}\dif s\\
&&+2\sqrt{\eps}\int^t_0\<v_\eps(s),\sigma(X^{\eps,h_\eps}(s))\dif W(s)\>_{\mR^m}\\
&&+\eps\int^t_0\|\sigma (X^{\eps,h_\eps}(s))\|^2_{L_2(l^2;\mR^m)}\dif s\\
&=:&I_1^\eps(t)+I_2^\eps(t)+I_3^\eps(t)+I_4^\eps(t)+I_5^\eps(t)+I_6^\eps(t).
\de

It is clear that by (\ref{Es8})
$$
I_1^\eps(t)\leq 0
$$
and
$$
I_2^\eps(t)\leq 2C_b\int^t_0|v_\eps(s)|^2\dif s.
$$
By BDG's inequality and {\bf (H2)} we also have
$$
\mE\left(\sup_{t\in[0,T]}|I_5^\eps(t)|\right)
+\mE\left(\sup_{t\in[0,T]}|I_6^\eps(t)|\right)\leq C\sqrt{\eps}.
$$
As estimating (\ref{Es4}) we have
$$
\mE\left(\sup_{s\in[0,t]}|I_3^\eps(s)|\right)
\leq \frac{1}{2}\mE\left(\sup_{s\in[0,t]}|v_\eps(s)|^2\right)
+C\int^t_0\mE|v_\eps(s)|^2\dif s.
$$

Set
$$
g(t):=\mE\left(\sup_{s\in[0,t]}|v_\eps(s)|^2\right).
$$
Then we have
$$
g(t)\leq \frac{1}{2}g(t)+C\sqrt{\eps}+\mE\left(\sup_{s\in[0,t]}|I_4^\eps(s)|\right)
+C\int^t_0\mE|v_\eps(s)|^2\dif s,
$$
which implies that
$$
g(t)\leq C\sqrt{\eps}+2\mE\left(\sup_{s\in[0,T]}|I_4^\eps(s)|\right)
+C\int^t_0g(s)\dif s.
$$
By Gronwall's inequality we get
\be
\mE\left(\sup_{s\in[0,T]}|v_\eps(s)|^2\right)\leq
C\sqrt{\eps}+C\mE\left(\sup_{s\in[0,T]}|I_4^\eps(s)|\right).\label{Es6}
\ee

We now deal with the hard term $I_4^\eps$. By It\^o's formula again, we have
\ce
\frac{1}{2}I_4^\eps(t)&=&\<v_\eps(t),w_\eps(t)\>_{\mR^m}
-\int^t_0w_\eps(s)\dif (K^{\eps,h_\eps}(s)-K^h(s))\\
&&-\int^t_0w_\eps(s)(b(X^{\eps,h_\eps}(s))-b(X^h(s)))\dif s\no\\
&&-\int^t_0w_\eps(s)(\sigma(X^{\eps,h_\eps}(s))\dot h_\eps(s)-\sigma(X^h(s))\dot h(s))\dif s\\
&&-\sqrt{\eps}\int^t_0w_\eps(s)\sigma(X^{\eps,h_\eps}(s))\dif W(s)\\
&=:&I_{41}^\eps(t)+I_{42}^\eps(t)+I_{43}^\eps(t)+I_{44}^\eps(t)+I_{45}^\eps(t).
\de
For any $\d>0$ and $R>0$, we have
\ce
P\left(\sup_{t\in[0,T]}|I_{41}^\eps(t)|\geq\d\right)
&=&P\left(\sup_{t\in[0,T]}|I_{41}^\eps(t)|\geq\d;
\sup_{t\in[0,T]}|v_\eps(t)|<R\right)\\
&&+P\left(\sup_{t\in[0,T]}|I_{41}^\eps(t)|\geq\d;
\sup_{t\in[0,T]}|v_\eps(t)|\geq R\right)\\
&\leq&P\left(\sup_{t\in[0,T]}|w_\eps(t)|\geq\d/R\right)
+P\left(\sup_{t\in[0,T]}|v_\eps(t)|\geq R\right).
\de
By Lemma \ref{Le3} and (\ref{Es5}) we know
$$
\lim_{\eps\to 0}
\mP\left(\sup_{t\in[0,T]}|I_{41}^\eps(t)|\geq\d\right)=0.
$$
Noting that
$$
\sup_{t\in[0,T]}\left|I_{42}^\eps(t)\right|\leq
\sup_{s\in[0,T]}|w_\eps(s)|\cdot (|K^{\eps,h_\eps}(\cdot)|^0_T+|K^h(\cdot)|^0_T),
$$
as above, we also have
$$
\sup_{t\in[0,T]}\left|I_{42}^\eps(t)\right|
\to 0\mbox{ in probability.}
$$
Similarly, we have
$$
\sup_{t\in[0,T]}\left|I_{43}^\eps(t)\right|
+\sup_{t\in[0,T]}\left|I_{44}^\eps(t)\right|\to 0\mbox{ in probability.}
$$
Moreover, by BDG's inequality we have
\ce
\sqrt{\eps}\mE\left(\sup_{t\in[0,T]}|I_{45}^\eps(t)|\right)
\leq C\sqrt{\eps}.
\de
Combining the above calculations, we get
$$
\sup_{t\in[0,T]}|I_4^\eps(t)|\to 0\mbox{ in probability.}
$$
It is easy to see by (\ref{Es5}) that
$$
\sup_{\eps\in(0,1)}\mE\left(\sup_{t\in[0,T]}|I_4^\eps(t)|^2\right)<+\infty.
$$
Hence
$$
\lim_{\eps\to 0}\mE\left(\sup_{t\in[0,T]}|I_4^\eps(t)|\right )=0.
$$
Substituting this into (\ref{Es6}) we obtain
\ce
\lim_{\eps\to 0}\mE\left(\sup_{s\in[0,T]}|v_\eps(s)|^2\right)=0.
\de
Thus, we have proven that for all $x\in \overline{D(A)}$
$$
\sup_{t\in[0,T]}|v_\eps(t,x)|^2\to 0,\ \ \ \mbox{ in probability}.
$$
We now strengthen it by Lemma \ref{Le5} to
$$
\xi_{n,\eps}:=\sup_{t\in[0,T],x\in D(A), |x|\leq n}|v_\eps(t,x)|^2\to 0,\ \ \ \mbox{ in probability}.
$$

Set
$$
\mD_n^\delta:=\overline{D(A)}\cap\{x\in\mR^m: |x|\leq n\}\cap \delta\mZ^m,
$$
where $\delta>0$ and $\delta\mZ^m$ denotes the grid in $\mR^m$ with edge length $\d$.
It is clear that there are only finite many points in $\mD_n^\delta$. Hence
$$
\xi^\delta_{n,\eps}:=\sup_{t\in[0,T], x\in\mD_n^\delta}|v_\eps(t,x)|^2\to 0,\ \ \ \mbox{ in probability}.
$$
For any $x\in \mR^m$, let $x_\delta$ denote the left-lower corner point in $\delta\mZ^m$ so that
$$
|x-x_\d|\leq \delta.
$$
Noting that
$$
\xi_{n,\eps}\leq 2\xi^\delta_{n,\eps}+2\sup_{x\in \overline{D(A)},
|x|\leq n}\sup_{t\in[0,T]}|v_\eps(t,x)-v_\eps(t,x_\d)|^2,
$$
we have for any $\beta>0$ and some $\a>0$
\ce
P(\xi_{n,\eps}>4\b)&\leq&P(\xi^\delta_{n,\eps}+\sup_{x\in \overline{D(A)},
|x|\leq n}\sup_{t\in[0,T]}|v_\eps(t,x)-v_\eps(t,x_\d)|^2\geq 2\beta)\\
&\leq& P(\xi^\delta_{n,\eps}>\b)
+P\left(\sup_{x\in \overline{D(A)}, |x|\leq n}\sup_{t\in[0,T]}|v_\eps(t,x)-v_\eps(t,x_\d)|^2>\b\right)\\
&\leq& P(\xi^\delta_{n,\eps}>\b)
+\mE\left(\sup_{x\in \overline{D(A)}, |x|\leq n}\sup_{t\in[0,T]}|v_\eps(t,x)-v_\eps(t,x_\d)|^2\right)/\b\\
&\leq& P(\xi^\delta_{n,\eps}>\b)+C\delta^\a/\b,
\de
where the last step is due to Lemma \ref{Le5} and Kolmogorov's criterion.
First letting $\delta$ be small enough, then $\eps$ go to zero, we then obtain
$$
\lim_{\eps\to 0}P(\xi_{n,\eps}>4\b)=0.
$$
which yields the desired result.
\end{proof}

\bl\label{Le6}
{\bf (LD)$_\mathbf{1}$} holds.
\el
\begin{proof}
Let $h_\eps$ be a sequence in $\cA^T_N$ converge to $h$ in distribution.
Since $\cD_N$ is compact and the law of $W$ is tight,
$\{h_\eps,W\}$ is tight in $\cD_N\times\mC_T(\mU)$ by the definition of tightness.
Without loss of generality, we assume the law of $\{h_\eps,W\}$ weakly converges to $\mu$.
Then the law of $h$ is just $\mu(\cdot,\mC_T(\mU))$. Indeed, for any bounded continuous function
$g$ on $\cD_N$, we have
\ce
\mE(g(h))=\lim_{n\rightarrow\infty}\mE(g(h_\eps))=\int_{\cD_N} g(h)\mu(\dif h,\mC_T(\mU)).
\de
By Skorohod's representation theorem, there are $(\tilde \Omega,\tilde P)$ and
$\{\tilde h_\eps,\tilde W^\eps\}$ and $\{\tilde h,\tilde W\}$ such that

(1) $(\tilde h_\eps,\tilde W^\eps)$ a.s. converges to $(\tilde h,\tilde W)$;

(2) $(\tilde h_\eps,\tilde W^\eps)$ has the same law as $(h_\eps,W)$;

(3) The law of $\{\tilde h,\tilde W\}$ is $\mu$, and the law of
$h$ is the same as $\tilde h$.

Using Lemma \ref{Le4}, we get
$$
\Phi(\frac{1}{\sqrt{\eps}}\int^\cdot_0\dot{\tilde h}_\eps\dif s+\tilde W^\eps)\to X^{\tilde h}, \ \ in\ probability,
$$
where $\Phi$ is the strong solution functional(cf. \cite{wa}). From this, we derive
$$
\Phi(\frac{1}{\sqrt{\eps}}\int^\cdot_0\dot{h}_\eps\dif s+W)\to X^{h}, \ \ in\ distribution.
$$
Thus, {\bf (LD)$_\mathbf{1}$} holds.
\end{proof}

Similar to the proof of Lemma \ref{Le4}, one can easily verify that
\bl\label{Le7}
{\bf (LD)$_\mathbf{2}$} holds.
\el

Thus, by Lemmas \ref{Le6}, \ref{Le7} and Theorem \ref{Th2}, we have proved Theorem \ref{Main}.

\end{document}